\documentclass[10pt]{article}
\usepackage{amsfonts,amsmath,amssymb,amsopn,amsthm,mdwlist}
\usepackage[cns,uns]{llconnex}

\input diagxy

\newcommand{\too}{\longrightarrow}

\newcommand{\oto}[1]{\stackrel{\underrightarrow{\stackrel{#1}{\,\,\,\,\hphantom{\too}}}}{}}  

\theoremstyle{definition}

\author{\sc Brian Day \footnote{Department of Mathematics, Macquarie University, NSW 2109, Australia.} }
\date{April 2, 2009}
\title{Middle-Four Maps and Net Categories}

\begin{document}

\maketitle

\begin{abstract}

We briefly relate the existence of a middle-four interchange map in a category with two monoidal structures, to the standard Cockett and Seely
notion of a weakly distributive category.

\end{abstract}

\vspace{1cm}
\section{Introduction}

\renewcommand{\d}{I_\parr}
\renewcommand{\e}{I_\tens}

From a middle-four map we shall obtain an example of a structure (on the category of $\d$-bimodules) with two weak distributivity maps:
\begin{align*}
	d^l &: A \tens (B \parr C) \too (A \tens B) \parr C \\
	d^r &: (B \parr C) \tens A \too B \parr (C \tens A)
\end{align*}
\noindent satisfying the eleven (planar) conditions (2)-(4) and (7)-(14) of \cite{CS}, but not necessarily \cite{CS} (1), (5), and (6). We call such a
structure a ``net'' category (see \cite{DP}) with the view that it is a proof-net category (\cite{DP}) without negation, or, equivalently,
a (planar) weakly distributive category \cite{CS} without the unit objects.

Moreover, we can talk about a middle-four category with negation, which satisfies four more suitable axioms (c.f. \cite{CS} \S4 (15)-(18)), so that,
on passing to the category of $\d$-bimodules, we get a proof-net category with a middle-four map.

\section{Middle-Four Categories}

First, we consider a category $V$ with two monoidal structures $(\tens, \e, a, l, r)$ and $(\parr, \d)$, where we have assumed that
the second structure is strictly monoidal. Suppose this category is equipped with a natural middle-four transformation:
\[ m : (A \parr B) \tens (C \parr D) \too (A \tens C) \parr (B \tens D) \]
and maps: \begin{align*}
	\mu &: \d \tens \d \too \d \\
	\eta &: \e \too \d
\end{align*}
satisfying the following four axioms:
\subsection*{M1}
\[ (a \parr a)m(m \tens 1) = m(1 \tens m)a : ((U \parr V) \tens (W \parr X)) \tens (Y \parr Z) \too (U \tens (W \tens Y)) \parr 
(V \tens (X \tens Z)) \]
\subsection*{M2}
\[ (m \parr 1)m = (1 \parr m)m \] from \[ (U \parr V \parr W) \tens (X \parr Y \parr Z) \] to \[ (U \tens X) \parr (V \tens Y) \parr (W \parr Z) \]
and, if we put $I = \e$ and $R = \d$, we want $(R,\mu,\eta)$ to be a $\tens$-monoid in $V$ such that both
\subsection*{M3}
\[ (1 \parr (\mu \tens 1))m(m \tens 1) = m \] from \[ ((A \parr R) \tens (B \parr R)) \tens (C \parr D) \] to \[ (( A \tens B) \tens C) \parr (R \tens D) \]
\subsection*{M4}
\[ ((1 \tens \mu) \parr 1)m(1 \tens m) = m \] from \[ (B \parr C) \tens ((R \parr A) \tens (R \parr B)) \] to
\[ (B \tens R) \parr ( C \tens (A \tens D)) \]

Moreover, a middle-four category with negation has maps $\gamma$ and $\tau$, where
\begin{align*}
	\gamma &: A \tens {A^\bot} \too R \\
	\tau &: I \too A^\bot \parr A
\end{align*}

satisfying four axioms similar to:
\[ \bfig

\node a(-500,+500)[(A \parr R) \tens I]
\node b(+500,+500)[(A \parr R) \tens (A^\bot \parr A)]
\node c(-500,0)[A]
\node d(+500,0)[(A \tens A^\bot) \parr (R \tens A)]
\node e(-500,-500)[R \tens A]
\node f(+500,-500)[R \parr (R \tens A)]

\arrow|a|[a`b;1 \tens \tau]
\arrow|l|[a`c;r]
\arrow|l|[c`e;\eta \tens 1]
\arrow|r|[b`d;m]
\arrow|r|[d`f;\gamma \parr 1]
\arrow|b|[f`e;=]

\efig \]
commutes (cf. \cite{CS} \S4).

\section{The Category of $R$-bimodules as a net category}

The category $V(R)$ of $R$-bimodules in $V$ has two induced (coherent) tensor products. These are:
\begin{itemize}
\item[(1)] $A \tens B$ with action: \begin{align*}
R \tens (A \tens B) \tens R &\oto{\sim} (R \tens A) \tens (B \tens R) \\
	&\oto{\sim} (R \tens A \tens I) \tens (I \tens B \tens R) \\
	&\oto{(1 \tens \eta) \tens (\eta \tens 1)} (R \tens A \tens R) \tens (R \tens B \tens R) \\
	&\oto{\alpha \tens \beta} A \tens B
\end{align*}
which is unitless, and
\item[(2)] $A \parr B$ with action 
\[ R \tens (A \parr B) \tens R \oto{\mbox{\tiny``$m$''}} (R \tens A \tens R) \parr (R \tens B \tens R) \oto{\alpha \parr \beta} A \parr B \]
which has the unit $R$.
\end{itemize}

Furthermore, it is easily checked that our basic transformation \[ m : (A \parr B) \tens (C \parr D) \too (A \tens C) \parr (B \tens D) \]
becomes an $R$-bimodule homomorphism. Thus, on putting $B = R$ and $C = R$ respectively, and using the bimodule actions, we get
two natural transformations:
\begin{align*}
	d^l &: A \tens (B \parr C) \too (A \tens B) \parr C \\
	d^r &: (B \parr C) \tens A \too B \parr (C \tens A)
\end{align*}
of $R$-bimodules.

We now check some \cite{CS} axioms on $d^l$ and $d^r$.

Axioms \cite{CS} (7) and (8) are trivial.

Axiom \cite{CS} (9):

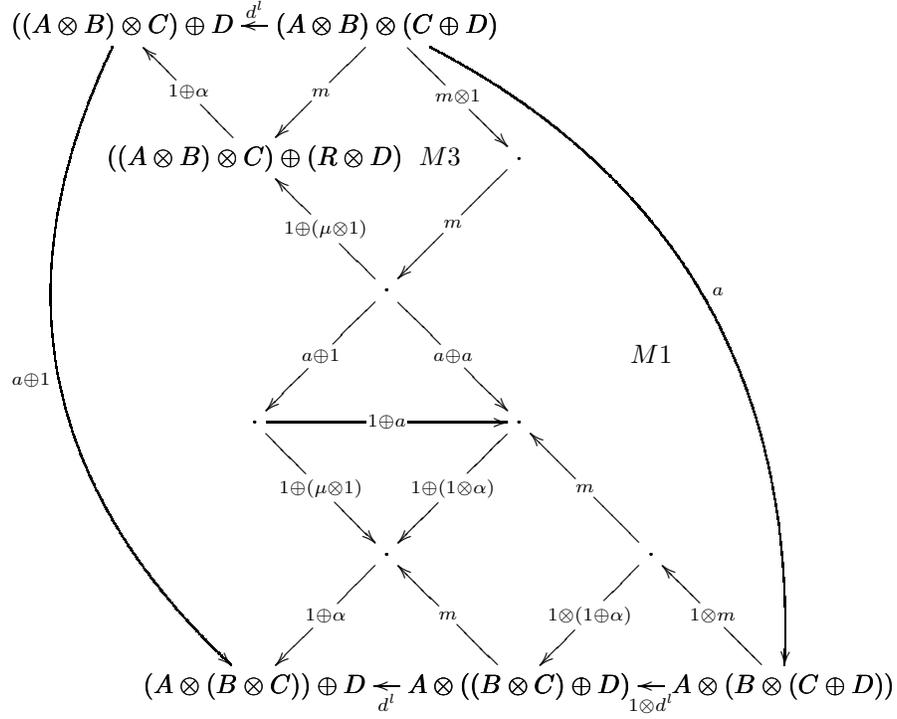
\[ \bfig

\node a(0,+1500)[(A \tens B) \tens (C \parr D)]
\node b(+1500,-1000)[A \tens (B \tens (C \parr D))]
\node c(-1000,+1500)[((A \tens B) \tens C) \parr D]
\node d(-500,-1000)[(A \tens (B \tens C)) \parr D]
\node e(+500,-1000)[A \tens ((B \tens C) \parr D)]
\node f(+500,+1000)[\cdot]
\node g(0,+500)[\cdot]
\node h(-500,+1000)[((A \tens B) \tens C) \parr (R \tens D)]
\node i(-500,0)[\cdot]
\node j(+500,0)[\cdot]
\node k(0,-500)[\cdot]
\node l(+1000,-500)[\cdot]

\place(+1000,+250)[M1]
\place(+200,+1000)[M3]

\arrow|m|[g`i;a \parr 1]
\arrow|m|[g`j;a \parr a]
\arrow|m|[i`j;1 \parr a]
\arrow|m|[i`k;1 \parr (\mu \tens 1)]
\arrow|m|[j`k;1 \parr (1 \tens \alpha)]
\arrow|m|[a`h;m]
\arrow|m|[a`f;m \tens 1]
\arrow|m|[f`g;m]
\arrow|m|[g`h;1 \parr (\mu \tens 1)]
\arrow|a|[a`c;d^l]
\arrow|m|[h`c;1 \parr \alpha]
\arrow|m|[l`e;1 \tens (1 \parr \alpha)]
\arrow|m|[e`k;m]
\arrow|m|[l`j;m]
\arrow|m|[b`l;1 \tens m]
\arrow|b|[b`e;1 \tens d^l]
\arrow|m|[k`d;1 \parr \alpha]
\arrow|b|[e`d;d^l]
\arrow|l|/{@{>}@/_5em/}/[c`d;a \parr 1]
\arrow|r|/{@{>}@/^5em/}/[a`b;a]

\efig \]

where $\alpha(1 \tens \alpha)a = \alpha(\mu \tens 1)$ sine $D$ is an $R$-bimodule.

Axiom \cite{CS} (10) is similar, but uses M1 and M4.

Axiom \cite{CS} (11) uses M2 and the action (*) of $R$ on $A \parr B$.

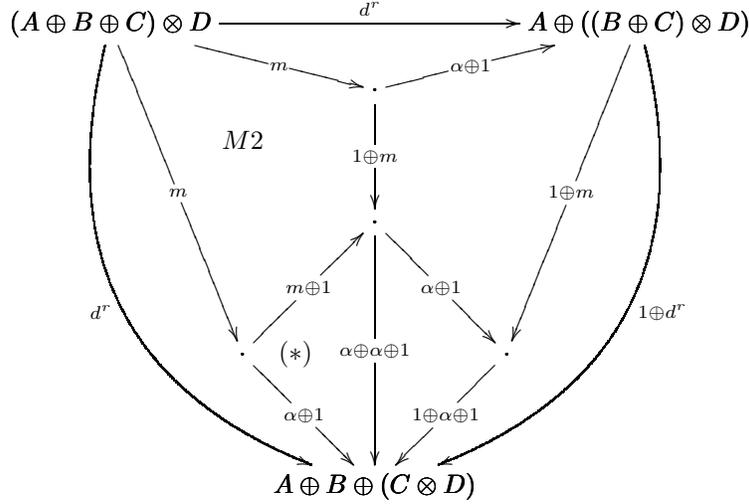
\[ \bfig

\node a(-1000,+750)[( A \parr B \parr C) \tens D]
\node b(+1000,+750)[A \parr (( B \parr C) \tens D)]
\node c(0,-1000)[A \parr B \parr (C \tens D)]
\node d(0,+500)[\cdot]
\node e(-500,-500)[\cdot]
\node f(0,0)[\cdot]
\node g(+500,-500)[\cdot]

\place(-500,+300)[M2]
\place(-300,-500)[(*)]

\arrow|a|[a`b;d^r]
\arrow|m|[a`d;m]
\arrow|m|[a`e;m]
\arrow|l|/{@{>}@/_5em/}/[a`c;d^r]
\arrow|m|[d`b;\alpha \parr 1]
\arrow|m|[b`g;1 \parr m]
\arrow|r|/{@{>}@/^5em/}/[b`c;1 \parr d^r]
\arrow|m|[g`c;1 \parr \alpha \parr 1]
\arrow|m|[f`c;\alpha \parr \alpha \parr 1]
\arrow|m|[e`c;\alpha \parr 1]
\arrow|m|[d`f;1 \parr m]
\arrow|m|[e`f;m \parr 1]
\arrow|m|[f`g;\alpha \parr 1]

\efig \]

Axiom \cite{CS} (12) is similar and still uses M2, while Axiom \cite{CS}(13) uses M2, and \cite{CS} (14) uses M1.

Thus, from any category with two monoidal structures related by a middle-four interchange transformation satisfying the standard axioms for such a 
map, we obtain an example of a net-category (that is, a proof-net category without negation) which is also equipped with a middle-four interchange map.

Communication on this article can be made through Micah McCurdy (Macquarie University, \verb$micah.mccurdy@gmail.com$) 
who kindly typed this manuscript. I would also like to thank
Ross Street for a helpful comment on the middle-four map.

\end{document}